\documentclass{article}
\usepackage{graphicx}
\usepackage{latexsym}
\usepackage{amssymb}
\usepackage{amsfonts}
\usepackage{amsmath}
\usepackage{indentfirst}
\usepackage{booktabs}
\usepackage[hang, small, bf]{caption}
\usepackage{float}
\usepackage{placeins}
\usepackage{latexsym}
\usepackage{tikz}

\newtheorem{proposition}{Proposition}[section]
\newtheorem{example}{Example}[section]

\newcommand{\cvd}{\hfill $\blacksquare$\bigskip}

\begin{document}
%
\title{Cross-bifix-free sets via Motzkin paths generation}
%
%
%

\date{}
\author{Elena~Barcucci\thanks{Dipartimento di Matematica e Informatica `U.Dini', Universit\`a degli Studi di Firenze, Viale
 G.B. Morgagni 65, 50134 Firenze, Italy (e-mail: elena.barcucci@unifi.it;
 stefano.bilotta@unifi.it; elisa.pergola@unifi.it; renzo.pinzani@unifi.it; jons10@hotmail.it).}
        \and Stefano~Bilotta$^*$
        \and Elisa~Pergola$^*$
        \and Renzo~Pinzani$^*$
        \and Jonathan~Succi$^*$}

\maketitle

\begin{abstract}
Cross-bifix-free sets are sets of words such that no prefix of any
word is a suffix of any other word. In this paper, we introduce a
general constructive method for the sets of cross-bifix-free
$q$-ary words of fixed length. It enables us to determine a
cross-bifix-free words subset which has the property to be
non-expandable.
\end{abstract}

\section{Introduction}
A \emph{cross-bifix-free set} of words (also called \emph{cross-bifix-free code}) is a set where, given any two words over an alphabet,
possibly the same, any prefix of the first one is not a suffix of the second one
and vice-versa.
Cross-bifix-free sets are involved in the study of frame synchronization which is an essential
requirement in a digital communication systems to establish and maintain a connection between a
transmitter and a receiver.

Analytical approaches to the synchronization acquisition process
and methods for the construction of sequences with the best
aperiodic autocorrelation properties \cite{1,2,3,4} have been the
subject of numerous analyses in the digital transmission.

The historical engineering approach started with the introduction
of bifix, a name proposed by J. L. Massey as acknowledged in
\cite{5}. It denotes a subsequence that is both a prefix and
suffix of a longer observed sequence.

In \cite{4} the notion of a \emph{distributed sequences} is introduced,
where the synchronization word is not a contiguous sequence of symbols
but is instead interleaved into the data stream. In \cite{6} is showed
that the distributed sequence entails a simultaneous search for a set of
synchronization words. Each word in the set of sequences is required to
be bifix-free. In addition, they arises a new requirement that no
prefix of any length of any word in the set is a suffix of any other word
in the set. This property of the set of synchronization words was termed
as \emph{cross-bifix-free}.

The problem of determining such sets is also related
to several other scientific applications, for instance in pattern
matching \cite{7} and automata theory \cite{8}.

Several methods for constructing cross-bifix-free sets have been recently proposed as in
\cite{9,10,11}. In particular, once the cardinality $q$ of the alphabet and the length $n$ of the
words are fixed, a matter is the construction of a cross-bifix-free set with the cardinality
as large as possible. An interesting method has been proposed in \cite{9} for words on
a binary alphabet. This specific construction reveals interesting connections to the Fibonacci
sequence of numbers. In a recent paper \cite{11} the authors revisit
the construction in \cite{9} and generalize it obtaining cross-bifix-free sets having greater cardinality
over an alphabet of any size $q$. They also show that their cross-bifix-free
sets have a cardinality close to the maximum possible. To our knowledge this is the best result in the literature about
the greatest size of cross-bifix-free sets.

For the sake of completeness we note that an intermediate step between
the original method \cite{9} and its generalization \cite{11} has been proposed in \cite{10} and it is constituted by
a different construction of binary cross-bifix-free sets based on lattice paths
which allows to obtain greater values of cardinality if compared to the ones in \cite{9}.

In this study, we revisit the construction in \cite{10}. We give a new construction of cross-bifix-free sets
that generalizes the construction of \cite{10} in order to extend the construction to $q$-ary alphabets for any $q$,
$q>2$. This approach enables us to obtain cross-bifix-free sets having greater
cardinality than the ones presented in \cite{11}, for the initial values of $n$.
This new result extends the theory of cross-bifix-free sets and it could be used to improve some technical applications.

This paper is organized as follows. In Section 2 we give some preliminaries and
describe the adopted notation. In Section 3 we present a new construction of cross-bifix-free sets in the $q$-ary
alphabet and in Section 4 we analyze the sizes of the sets of our construction in comparison to the ones in the literature.

\section{Basic definitions and notations}
Let $\mathbb{Z}_q=\{0,1,\cdots,q-1\}$ be an alphabet of $q$ elements. A (finite) sequence of
elements in $\mathbb{Z}_q$ is called (finite) \emph{word}. The set of all words over $\mathbb{Z}_q$
having length $n$ is denoted by $\mathbb{Z}^n_q$.
A consecutive sequence of $m$ element $a \in \mathbb{Z}_q$
is denoted by the short form $a^m$.
Let $w \in \mathbb{Z}_q^n$, then ${|w|}_a$ denotes the
number of occurrences of $a$ in $w$, being $a \in \mathbb{Z}_q$. Let
$w=uzv$ then $u$ is called a \emph{prefix} of $w$ and $v$
is called a \emph{suffix} of $w$. A \emph{bifix} of $w$ is
a subsequence of $w$ that is both its prefix and suffix.

A word $w \in \mathbb{Z}_q^n$ is said to be \emph{bifix-free} or
\emph{unbordered} \cite{12} if and only if no prefix of
$w$ is also a suffix of $w$. Therefore, $w$ is
bifix-free if and only if $w \neq uzu$, being $u$ any
necessarily non-empty word and $z$ any word. Obviously, a
necessary condition for $w$ to be bifix-free is that the
first and the last letters of $w$ must be different.
\begin{example}
In $\mathbb{Z}_2=\{0,1\}$, the word $111010100$ of length $n=9$ is
bifix-free, while the word $101001010$ contains two bifixes, $10$
and $1010$.
\end{example}

Let $BF_q(n)$ denote the set of all bifix-free words of length $n$
over an alphabet of fixed size $q$ (for more details about this topic see \cite{12}).

Given $q>1$ and $n>1$, two distinct words $w,w'
\in BF_q(n)$ are said to be \emph{cross-bifix-free} \cite{6} if
and only if no strict prefix of $w$ is also a suffix of
$w'$ and vice-versa.

\begin{example}\label{esemp}
The binary words $111010100$ and $110101010$ in $BF_2(9)$ are
cross-bifix-free, while the binary words $111001100$ and
$110011010$ in $BF_2(9)$ have the cross-bifix $1100$.
\end{example}

A subset of $BF_q(n)$ is said to be a \emph{cross-bifix-free set} if
and only if for each $w, w'$, with $w \neq
w'$, in this set, $w$ and $w'$ are
cross-bifix-free. This set is said to be \emph{non-expandable} on
$BF_q(n)$ if and only if the set obtained by adding any other word in $BF_q(n)$
is not a cross-bifix-free set. A non-expandable cross-bifix-free
set on $BF_q(n)$ having maximal cardinality is called a
\emph{maximal cross-bifix-free set} on $BF_q(n)$.

\medskip

In a recent paper \cite{11} the authors provide a general construction of cross-bifix-free sets
over a $q$-ary alphabet. Below, we recall such generation for the family of cross-bifix-free sets in $\mathbb{Z}_q^n$.

For any $2 \le k \le n-2$, the cross-bifix-free set $\mathcal S_{k,q}(n)$ in \cite{11}
is the set of all words $s=s_1 s_2 \cdots s_n$ in $\mathbb{Z}_q^n$ that satisfy the following two properties:
\begin{itemize}
\item[1)] $s_1 = \dots = s_k = 0$, $s_{k+1} \ne 0$ and $s_n \ne 0$,
\item[2)] the subsequence $s_{k+2} \dots s_{n-1}$ does not contain $k$ consecutive 0's.
\end{itemize}

\medskip

Let

\begin{equation*}
F_{k,q}(n)=
\begin{cases}
q^n & \text{if $0 \le n < k$,}\\
(q-1)\sum_{l=1}^k F_{k,q}(n-l) & \text{if $n \ge k$},
\end{cases}
\end{equation*}

\smallskip
\noindent
be the sequence enumerating the words in $\mathbb{Z}_q^n$ avoiding $k$ consecutive zero's
\cite{13}. Then, from the above definition of $\mathcal S_{k,q}(n)$, we have

\[|S_{n,q}^{(k)}|=(q-1)^2 F_{k,q}(n-k-2)\ .\]

For any fixed $n$ and $q$, the largest size of $|S_{n,q}^{(k)}|$ is denoted by $S(n,q)$ and it is given by the following expression as in \cite{11}
\[ S(n,q) = \max\{(q-1)^2 F_{k,q}(n-k-2): 2 \leq k \leq n-2 \}.\]

This result allows to obtain non-expandable cross-bifix-free
sets in the $q$-ary alphabet having cardinality close to the maximum.

\medskip

In the present paper we introduce an alternative constructive method for the generation of cross-bifix-free set in $\mathbb{Z}_q$. Our approach is based on the study of lattice path in the discrete plane and it moves from the construction in \cite{10}.

Each word $w \in \mathbb{Z}_q^n$ can be represented as a lattice path of $\mathbb{N}^2$ running from $(0,0)$ to $(n,0)$ having the following properties:
\begin{itemize}
\item[-] the element $0$ corresponds to a \emph{fall step} which is defined by $(1,-1)$,
\item[-] the element $1$ corresponds to a \emph{rise step} which is defined by $(1,1)$,
\item[-] the elements $2,\dots,q-1$ correspond respectively to a \emph{colored level step} which is defined by $(1,0)$ and it is labeled by one of the $q-2$ fixed colors.
\end{itemize}

For example, in Table \ref{tab1} and Table \ref{tab2} is showed an equivalence between elements
and steps of lattice paths in the alphabets $\mathbb{Z}_3$ and $\mathbb{Z}_4$, respectively.

\begin{table}[H]
\caption{Equivalence between symbols and steps for \\ $\mathbb{Z}_3 = \{0,1,2\}$.}
\centering
\begin{tabular}{c c c c}
\toprule
Symbol & Step & Color & Representation\\
\midrule
$0$ & $(1,-1)$ & - & \begin{tikzpicture}[scale=0.5] \draw[help lines, dotted] (0,0) grid (1,1); \draw [very thick] (0,1) -- (1,0); \end{tikzpicture}\\
$1$ & $(1,1)$ & - & \begin{tikzpicture}[scale=0.5] \draw[help lines, dotted] (0,0) grid (1,1); \draw [very thick] (0,0) -- (1,1); \end{tikzpicture}\\
$2$ & $(1,0)$ & Black & \begin{tikzpicture}[scale=0.5] \draw[help lines, dotted] (0,0) grid (1,1); \draw [very thick] (0,0) -- (1,0); \end{tikzpicture}\\
\bottomrule
\end{tabular}\label{tab1}
\end{table}

\begin{table}[H]
\caption{Equivalence between symbols and steps for \\ $\mathbb{Z}_4 = \{0,1,2,3\}$.}
\centering
\begin{tabular}{c c c c}
\toprule
Symbol & Step & Color & Representation\\
\midrule
$0$ & $(1,-1)$ & - & \begin{tikzpicture}[scale=0.5] \draw[help lines, dotted] (0,0) grid (1,1); \draw [very thick] (0,1) -- (1,0); \end{tikzpicture}\\
$1$ & $(1,1)$ & - & \begin{tikzpicture}[scale=0.5] \draw[help lines, dotted] (0,0) grid (1,1); \draw [very thick] (0,0) -- (1,1); \end{tikzpicture}\\
$2$ & $(1,0)$ & Red & \begin{tikzpicture}[scale=0.5] \draw[help lines, dotted] (0,0) grid (1,1); \draw [red, very thick] (0,0) -- (1,0); \end{tikzpicture}\\
$3$ & $(1,0)$ & Green & \begin{tikzpicture}[scale=0.5] \draw[help lines, dotted] (0,0) grid (1,1); \draw [green, very thick] (0,0) -- (1,0); \end{tikzpicture}\\
\bottomrule
\end{tabular}\label{tab2}
\end{table}

From now on, we will refer interchangeably to words or their
graphical representations on the discrete plane, that is paths. The definition of bifix-free and cross-bifix-free can be easily
extended to paths.

\medskip

A \emph{$k$-colored Motzkin path} of length $n$ is a lattice path of $\mathbb{N}^2$ running from $(0,0)$ to $(n,0)$ that never goes below the $x$-axis and whose admitted steps are rise steps, fall steps and $k$-colored level steps (for more details about this copy see \cite{14}).

\begin{figure}
\centering
\caption{Words $121002$, $100212$ and the equivalent paths. The first one is a Motzkin word.}
\begin{tikzpicture}[scale=0.5]
\draw (0,0) -- (1,1) -- (2,1) -- (3,2) -- (4,1) -- (5,0) -- (6,0);
\node [below] at (0.5,-1) {$1$};
\node [below] at (1.5,-1) {$2$};
\node [below] at (2.5,-1) {$1$};
\node [below] at (3.5,-1) {$0$};
\node [below] at (4.5,-1) {$0$};
\node [below] at (5.5,-1) {$2$};
\end{tikzpicture}
\quad
\begin{tikzpicture}[scale=0.5]
\draw (0,0) -- (1,1) -- (2,0) -- (3,-1) -- (4,-1) -- (5,0) -- (6,0);
\node [below] at (0.5,-1) {$1$};
\node [below] at (1.5,-1) {$0$};
\node [below] at (2.5,-1) {$0$};
\node [below] at (3.5,-1) {$2$};
\node [below] at (4.5,-1) {$1$};
\node [below] at (5.5,-1) {$2$};
\end{tikzpicture}\label{f1}
\end{figure}
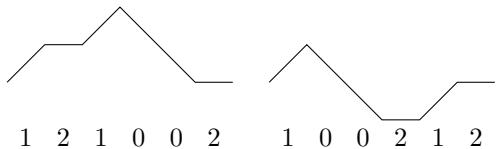

For example, the left side of Fig. \ref{f1} shows a Motzkin path in $\mathbb{Z}_3$ having length 6, while the path in its right side is not a Motzkin path since it crosses the $x$-axis.

We denote by $\mathcal{M}_k(n)$ the set of all $k$-colored Motzkin paths of length $n$, and let $M_k(n)$ be the size of $\mathcal{M}_k(n)$.

\medskip

\begin{proposition}
For any $n \ge 0$ and $k \ge 3$, $M_k(n)$ is given by the following expression
\begin{equation*}
\begin{split}
M_k(n+1) & = k M_k(n) + \sum_{i=0}^{n-1}M_k(i) M_k(n-1-i)
\end{split}
\end{equation*}
with $M_k(0)  = 1$ and
$M_k(1) = k$.
\end{proposition}

\smallskip

\emph{Proof.} \quad If $n = 0$, $\mathcal{M}_k(n)$ contains the empty path only, then $M_k(0) = 1$. If $n = 1$, $\mathcal{M}_k(n)$ only contains those paths obtained by a level step, thus $M_k(1) = k$.\\
Let $n \ge 1$ and $w \in \mathcal{M}_k(n+1)$. There are two cases: $w$ begins with a level step or $w$ begins with a rise step. In the first case we have that $w = h \alpha$ where $h$ is a level step and $\alpha \in \mathcal{M}_k(n)$, then the number of this first kind of paths is equal to $k M_k(n)$.

Otherwise, we have that $w = u  \alpha  d \beta$
where $u$ is a rise step, $d$ is a fall step, $\alpha \in \mathcal{M}_k(i)$ and $\beta \in \mathcal{M}(n-1-i)$ with $0 \le i \le n-1$. Then the number of this latter kind of paths is equal to $\sum_{i=0}^{n-1}M_k(i) M_k(n-1-i)$.

Thus,
\[M_k(n+1) = k M_k(n) + \sum_{i=0}^{n-1}M_k(i) M_k(n-1-i).\]
\cvd

A word $w \in \mathbb{Z}_q^n$ is called \emph{$(q-2)$-colored Motzkin word} if the equivalent lattice path is a $(q-2)$-colored Motzkin path.

\smallskip

For our purposes, it is useful to denote by $\hat{\mathcal{M}}_{q-2}(n)$ the set of all \emph{elevated $(q-2)$-colored Motzkin words} of length $n$, defined as

\[ \hat{\mathcal{M}}_{q-2}(n) = \{1 \alpha 0 : \alpha \in \mathcal{M}_{q-2}(n - 2) \}.\]

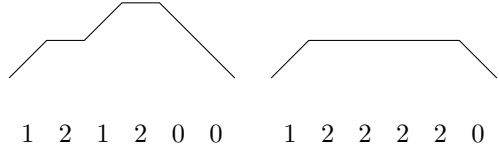
\begin{figure}
\centering
\caption{An example of elevated Motzkin words}
\begin{tikzpicture}[scale=0.5]
\draw (0,0) -- (1,1) -- (2,1) -- (3,2) -- (4,2) -- (5,1) -- (6,0);
\node [below] at (0.5,-1) {$1$};
\node [below] at (1.5,-1) {$2$};
\node [below] at (2.5,-1) {$1$};
\node [below] at (3.5,-1) {$2$};
\node [below] at (4.5,-1) {$0$};
\node [below] at (5.5,-1) {$0$};
\end{tikzpicture}
\quad
\begin{tikzpicture}[scale=0.5]
\draw (0,0) -- (1,1) -- (2,1) -- (3,1) -- (4,1) -- (5,1) -- (6,0);
\node [below] at (0.5,-1) {$1$};
\node [below] at (1.5,-1) {$2$};
\node [below] at (2.5,-1) {$2$};
\node [below] at (3.5,-1) {$2$};
\node [below] at (4.5,-1) {$2$};
\node [below] at (5.5,-1) {$0$};
\end{tikzpicture}\label{f2}
\end{figure}

For example, in Fig. \ref{f2} two words in $\hat{\mathcal{M}}_1(6)$ are depicted.

In the next section of the present paper we are interested in
determining one among all the possible non-expandable
cross-bifix-free sets of words of fixed length $n>1$ on $\mathbb{Z}_q^n$. We denote this set by $\mathcal{CBFS}_q(n)$.

\section{On the non-expandability of $\mathcal{CBFS}_q(n)$}
In this section we define the set $\mathcal{CBFS}_q(n)$ which is formed by the union of three sets of $(q-2)$-colored Motzkin paths denoted by
$\mathcal{A}_q(n), \mathcal{B}_q(n)$ and $\mathcal{C}_q(n)$, with $q \ge 3$ and $n \ge 3$, respectively.

\medskip

\begin{figure}
\caption{Graphical representation of the set $\mathcal{A}_q(n)$, $n \ge 3$}
\centering
\begin{tikzpicture}[scale=0.4]
\draw (0,0) to [out=45, in=180] (2,2) to [out=0, in=135] (3,1) to [out=45, in=180] (5,3) to [out=0, in=135] (8,0) -- (0,0);
\draw (8,0) -- (9,1);
\draw (9,1) to [out=45, in=180] (11,3) to [out=0, in=135] (12,2) to [out=45, in=180] (14,4) to [out=0, in=135] (17,1) -- (9,1);
\draw (17,1) -- (18,0);
\node [below] at (4,-1) {$\alpha \in \mathcal{M}_{q-2}(i)$};
\node [below] at (13,-1) {$\beta \in \hat{\mathcal{M}}_{q-2}(n-i)$};
\end{tikzpicture}\label{A}
\end{figure}
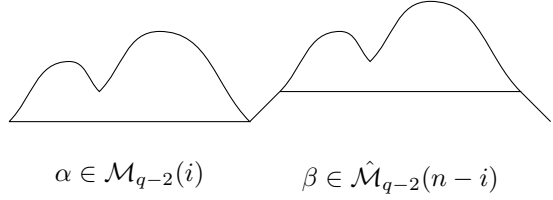

Let
\begin{equation*}
\begin{split}
\mathcal{A}_q(n) = & \left\{ \alpha \beta : \alpha \in \mathcal{M}_{q-2}(i), \beta \in \hat{\mathcal{M}}_{q-2}(n-i) \right\} \setminus
 \left\{ \alpha \beta : \alpha, \beta \in \hat{\mathcal{M}}_{q-2}\left(\frac{n}{2}\right) \right\}
\end{split}
\end{equation*}
with $0 \le i \le \left\lfloor \frac{n}{2} \right\rfloor$, be the set of words composed by a $(q-2)$-colored Motzkin word $\alpha$ of length $i$, and a elevated $(q-2)$-colored Motzkin word $\beta$ of length $n-i$ (see Fig. \ref{A}). If $n$ is even, we need to remove the words composed by two elevated subwords of the same length. On the other side, if $n$ is odd, we assume the set $\left\{ \alpha \beta : \alpha, \beta \in \hat{\mathcal{M}}_{q-2}\left(\frac{n}{2}\right) \right\}$ empty, since it does not exists any path of non-integer length.

Then, the enumeration of the set $\mathcal{A}_q(n)$ is given by the following expression
\begin{equation*}
\begin{split}
|\mathcal{A}_q(n)| = & \sum_{i=0}^{\left\lfloor n/2 \right\rfloor} M_{q-2}(i) M_{q-2}(n-i-2) - \left[M_{q-2}\left(\frac{n}{2} - 2\right)\right]^2.
\end{split}
\end{equation*}

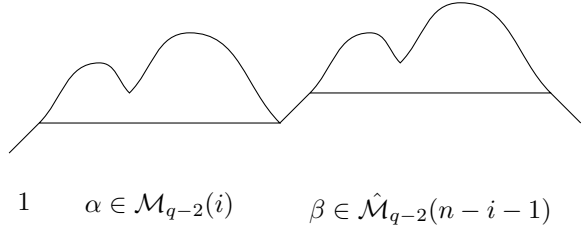
\begin{figure}
\caption{Graphical representation of the set $\mathcal{B}_q(n)$, $n \ge 3$}
\centering
\begin{tikzpicture}[scale=0.4]
\draw (0,0) -- (1,1);
\draw (1,1) to [out=45, in=180] (3,3) to [out=0, in=135] (4,2) to [out=45, in=180] (6,4) to [out=0, in=135] (9,1) -- (1,1);
\draw (9,1) -- (10,2);
\draw (10,2) to [out=45, in=180] (12,4) to [out=0, in=135] (13,3) to [out=45, in=180] (15,5) to [out=0, in=135] (18,2) -- (10,2);
\draw (18,2) -- (19,1);
\node [below] at (0.5,-1) {$1$};
\node [below] at (5,-1) {$\alpha \in \mathcal{M}_{q-2}(i)$};
\node [below] at (14,-1) {$\beta \in \hat{\mathcal{M}}_{q-2}(n-i-1)$};
\end{tikzpicture}\label{B}
\end{figure}

Let
\begin{equation*}
\mathcal{B}_q(n) = \left\{ 1 \alpha \beta : \alpha \in \mathcal{M}_{q-2}(i), \beta \in \hat{\mathcal{M}}_{q-2}(n-i-1) \right\}
\end{equation*}
with $0 \le i \le \left\lfloor \frac{n}{2} \right\rfloor - 1$, be the set of words composed by a rise step, a $(q-2)$-colored Motzkin word $\alpha$ of length $i$, and a elevated $(q-2)$-colored Motzkin word $\beta$ of length $n-i-1$ (see Fig. \ref{B}).

Then, the enumeration of the set $\mathcal{B}_q(n)$ is given by the following expression
\begin{equation*}
|\mathcal{B}_q(n)| = \sum_{i=0}^{\left\lfloor n/2 \right\rfloor - 1} M_{q-2}(i) M_{q-2}(n-i-3).
\end{equation*}

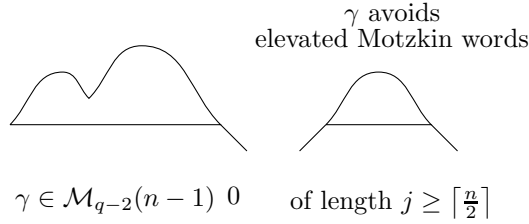
\begin{figure}
\caption{Graphical representation of the set $\mathcal{C}_q(n)$, $n \ge 3$}
\centering
\begin{tikzpicture}[scale=0.35]
\draw (0,1) to [out=45, in=180] (2,3) to [out=0, in=135] (3,2) to [out=45, in=180] (5,4) to [out=0, in=135] (8,1) -- (0,1);
\draw (8,1) -- (9,0);
\draw (11,0) -- (12,1);
\draw (12,1) to [out=45, in=180] (14,3) to [out=0, in=135] (16,1) -- (12,1);
\draw (16,1) -- (17,0);
\node [below] at (4,-1) {$\gamma \in \mathcal{M}_{q-2}(n-1)$};
\node [below] at (8.5,-1) {$0$};
\node [below] at (14.5,6) {$\gamma$ avoids};
\node [below] at (14.5,5) {elevated Motzkin words};
\node [below] at (14.5,-1) {of length $j \ge \left\lceil \frac{n}{2} \right\rceil$};
\end{tikzpicture}\label{C}
\end{figure}

Let
\begin{equation*}
\mathcal{C}_q(n) =\left\{ \gamma 0 : \gamma \in \mathcal{M}_{q-2}(n-1), \gamma \ne u \beta v, \beta \in \hat{\mathcal{M}}_{q-2}(j) \right\}
\end{equation*}
with $j \ge \left\lceil \frac{n}{2} \right\rceil$, be the set of words composed by a $(q-2)$-colored Motzkin word $\gamma$ of length $n-1$ that avoids elevated $(q-2)$-colored Motzkin words of length $j$, and a fall step (see Fig. \ref{C}).

Then, the enumeration of the set $\mathcal{C}_q(n)$ is given by the following expression
\begin{equation*}
\begin{split}
&|\mathcal{C}_q(n)| =  M_{q-2}(n-1) -
 \sum_{k = \left\lceil n/2 \right\rceil}^{n-1} \sum_{i=0}^{n-1-k} M_{q-2}(i) M_{q-2}(k-2) M_{q-2}(n-1-i-k).
\end{split}
\end{equation*}

Note that, in order to obtain the size $|\mathcal{C}_q(n)|$ we need to subtract from all words $\gamma$ of length $n-1$ those containing a elevated Motzkin subword $\beta$ of length greater than or equal to $\left\lceil n/2 \right\rceil$, and $\gamma$ can contain one of those subwords at most. Then, for $k = \left\lceil n/2 \right\rceil,\dots, n-1$ we need to remove the words $u \beta v$, with $u \in \mathcal{M}_{q-2}(i)$, $\beta \in \hat{\mathcal{M}}_{q-2}(k)$, $v \in \mathcal{M}_{q-2}(n-1-i-k)$ and $0 \le i \le n-1-k$.

\medskip

\begin{figure}[h!]
\caption{Graphical representation of the set $\mathcal{CBFS}_3(4)$}
\centering
\begin{tikzpicture}[scale=0.3]
\draw [out=45, in=180] (0,0) -- (1,1) -- (2,1) -- (3,1) -- (4,0);
\node [below] at (0.5,-1) {$1$};
\node [below] at (1.5,-1) {$2$};
\node [below] at (2.5,-1) {$2$};
\node [below] at (3.5,-1) {$0$};
\end{tikzpicture}
\quad
\begin{tikzpicture}[scale=0.3]
\draw [out=45, in=180] (0,0) -- (1,1) -- (2,2) -- (3,1) -- (4,0);
\node [below] at (0.5,-1) {$1$};
\node [below] at (1.5,-1) {$1$};
\node [below] at (2.5,-1) {$0$};
\node [below] at (3.5,-1) {$0$};
\end{tikzpicture}
\quad
\begin{tikzpicture}[scale=0.3]
\draw [out=45, in=180] (0,0) -- (1,0) -- (2,1) -- (3,1) -- (4,0);
\node [below] at (0.5,-1) {$2$};
\node [below] at (1.5,-1) {$1$};
\node [below] at (2.5,-1) {$2$};
\node [below] at (3.5,-1) {$0$};
\end{tikzpicture}
\quad
\begin{tikzpicture}[scale=0.3]
\draw [out=45, in=180] (0,0) -- (1,0) -- (2,0) -- (3,1) -- (4,0);
\node [below] at (0.5,-1) {$2$};
\node [below] at (1.5,-1) {$2$};
\node [below] at (2.5,-1) {$1$};
\node [below] at (3.5,-1) {$0$};
\end{tikzpicture}
\\
\quad
\\
\begin{tikzpicture}[scale=0.3]
\draw [out=45, in=180] (0,0) -- (1,1) -- (2,2) -- (3,2) -- (4,1);
\node [below] at (0.5,-1) {$1$};
\node [below] at (1.5,-1) {$1$};
\node [below] at (2.5,-1) {$2$};
\node [below] at (3.5,-1) {$0$};
\end{tikzpicture}
\quad
\begin{tikzpicture}[scale=0.3]
\draw [out=45, in=180] (0,0) -- (1,1) -- (2,1) -- (3,2) -- (4,1);
\node [below] at (0.5,-1) {$1$};
\node [below] at (1.5,-1) {$2$};
\node [below] at (2.5,-1) {$1$};
\node [below] at (3.5,-1) {$0$};
\end{tikzpicture}
\quad
\begin{tikzpicture}[scale=0.3]
\draw [out=45, in=180] (0,1) -- (1,1) -- (2,1) -- (3,1) -- (4,0);
\node [below] at (0.5,-1) {$2$};
\node [below] at (1.5,-1) {$2$};
\node [below] at (2.5,-1) {$2$};
\node [below] at (3.5,-1) {$0$};
\end{tikzpicture}\label{cbfs}
\end{figure}
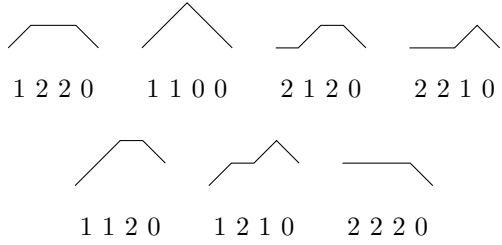

At this point, we define the set $\mathcal{CBFS}_q(n)$ as follows
\begin{equation*}
\mathcal{CBFS}_q(n) = \mathcal{A}_q(n) \cup \mathcal{B}_q(n) \cup \mathcal{C}_q(n)
\end{equation*}
that is the union of the above described sets. For instance, in Fig. \ref{cbfs} the set $\mathcal{CBFS}_3(4)$ is depicted.

\medskip

\begin{proposition}
The set $\mathcal{CBFS}_q(n)$ is a cross-bifix-free set on $BF_q(n)$, for any $q \ge 3$ and $n \ge 3$.
\end{proposition}

\smallskip

\emph{Proof.} \quad Let $w, w' \in \mathcal{CBFS}_q(n)$. Let $u$ be a prefix of $w$, and $v$ be a suffix of $w'$ such that $|u| = |v|$. We need to check that in each case the prefix $u$ does not match with the suffix $v$.
\begin{enumerate}
\item Let $w \in \mathcal{A}_q(n)$ and $w' \in \mathcal{A}_q(n) \cup \mathcal{B}_q(n)$.\\
For each prefix $u$ of $w$ we have $|u|_0 \le |u|_1$ and if $|u| > \lfloor \frac{n}{2} \rfloor$, then $|u|_0 < |u|_1$. For each suffix $v$ of $w'$ we have $|v|_0 \ge |v|_1$ and if $|v| < \lfloor \frac{n+1}{2} \rfloor$, then $|v|_0 > |v|_1$.\\
Let $|u| = |v| = l$, if either $l < \lfloor \frac{n+1}{2} \rfloor$ or $l > \lfloor \frac{n}{2} \rfloor$, then $u$ does not match with $v$. So we have to check the case $\lfloor \frac{n+1}{2} \rfloor \le l \le \lfloor \frac{n}{2} \rfloor$.

If $n$ is odd, it does not exist an integer $l$ satisfying $\lfloor \frac{n+1}{2} \rfloor \le l \le \lfloor \frac{n}{2} \rfloor$, otherwise if $n$ is even, the case $\lfloor \frac{n+1}{2} \rfloor \le l \le \lfloor \frac{n}{2} \rfloor$ is verified only for $l = \frac{n}{2}$. Therefore let $n$ be even and $l = \frac{n}{2}$. In this case $|u|_0 \le |u|_1$ and $|v|_0 \ge |v|_1$. At this point $u$ can match with $v$ only if $|v|_0 = |v|_1$, and this can happen only if $v$ is a elevated Motzkin word. Suppose now that $u = v$, so $u$ should be a elevated Motzkin word too, and they have both length $\frac{n}{2}$. In this case, $w$ should be a word composed of two elevated Motzkin subwords of the same length, but such a word does not exists in $\mathcal{CBFS}_q(n)$ since the set $\left\{ \alpha \beta : \alpha, \beta \in \hat{\mathcal{M}}_{q-2}\left(\frac{n}{2}\right) \right\}$ is not included in it, thus $u$ does not match with $v$.\\

\item Let $w \in \mathcal{B}_q(n)$ and $w' \in \mathcal{A}_q(n) \cup \mathcal{B}_q(n)$.\\
For each prefix $u$ of $w$ we have $|u|_0 < |u|_1$, and for each suffix $v$ of $w'$ we have $|v|_0 \ge |v|_1$, thus $u$ does not match with $v$.\\

\item Let $w \in \mathcal{C}_q(n)$ and $w' \in \mathcal{A}_q(n) \cup \mathcal{B}_q(n)$.\\
For each prefix $u$ of $w$ we have $|u|_0 \le |u|_1$. For each suffix $v$ of $w'$ we have $|v|_0 \ge |v|_1$ and if $|v| < \lfloor \frac{n+1}{2} \rfloor$, then $|v|_0 > |v|_1$.\\
Let $|u| = |v| = l$. If $l < \lfloor \frac{n+1}{2} \rfloor$, then $u$ does not match with $v$. So we have to check the case
$ l \ge \lfloor \frac{n+1}{2} \rfloor$.
In this case $v$ contains a elevated Motzkin subword of length $\lfloor \frac{n+1}{2} \rfloor = \lceil \frac{n}{2} \rceil$ at least, and $u$ does not match with $v$, since $u$ avoids such subwords.\\

\item Let $w \in \mathcal{CBFS}_q(n)$ and $w' \in \mathcal{C}_q(n)$.\\
For each prefix $u$ of $w$ we have $|u|_0 \le |u|_1$, and for each suffix $v$ of $w'$ we have $|v|_0 > |v|_1$, thus $u$ cannot match with $v$.\\

\end{enumerate}
We proved that $\mathcal{CBFS}_q(n)$ is a cross-bifix-free set on $BF_q(n)$, for any $q \ge 3$ and $n \ge 3$.
\cvd

\medskip

\begin{proposition}
The set $\mathcal{CBFS}_q(n)$ is a non-expandable cross-bifix-free set on $BF_q(n)$, for any $q \ge 3$ and $n \ge 3$.
\end{proposition}

\smallskip

\emph{Proof.} \quad Let $w \in BF_q(n) \setminus \mathcal{CBFS}_q(n)$ and $W = \mathcal{CBFS}_q(n) \cup \{w\}$.
If $w$ begins with $0$ then $W$ is not cross-bifix-free since any word in $\mathcal{CBFS}_q(n)$ ends with $0$. If $w$ ends with $1$ then $W$ is not cross-bifix-free since any word in $\mathcal{A}_q(n)$ begins with $1$. If $w$ ends with a letter $k \ne 0,1$ then $W$ is not cross-bifix-free since the suffix $k$ of $w$ matches, for instance, with the prefix $k$ of the word $k^{n-1}0 \in \mathcal{C}_q(n)$. Consequently we have to consider $w$ as a word beginning with a non-zero letter and ending with $0$.\\

Let $h = |w|_1 - |w|_0$ be the ordinate of the last point of the path corresponding to $w$. We now need to distinguish three different cases: $h > 0$, $h < 0$ and $h = 0$.\\

If $h > 0$, $w$ can be written as (see Fig. \ref{fig:hpos})
\[ w = \phi \; 1 \; \mu_1 \; 1 \; \mu_2 \; \cdots \; 1 \; \mu_h, \]
where $\phi$ is a word satisfying $|\phi|_1 = |\phi|_0$ and not beginning with $0$, and $\mu_1, \dots, \mu_h$ are $(q-2)$-colored Motzkin words with $\mu_h$ non-empty as $w$ ends with $0$.

\begin{figure}[h!]
\caption{Graphical representation of $w$, in the case $h > 0$}
\label{fig:hpos}
\centering
\begin{tikzpicture}[scale=0.3]
\draw (0,0) to [out=45, in=180] (2,2) to [out=0, in=180] (6,-2) to [out=0, in=225] (8,0) -- (0,0);
\draw (8,0) -- (9,1);
\draw (9,1) to [out=45, in=180] (11,3) to [out=0, in=135] (13,1) -- (9,1);
\draw (13,1) -- (14,2);
\draw (14,2) to [out=45, in=180] (16,4) to [out=0, in=135] (18,2) -- (14,2);
\draw [dashed] (18,2) -- (20,4);
\draw (20,4) -- (21,5);
\draw (21,5) to [out=45, in=180] (23,7) to [out=0, in=135] (25,5) -- (21,5);
\draw [dashed] (0,0) -- (25,0) -- (25,5);
\node [above] at (4,-4) {$\phi$};
\node [above] at (8.5,-4) {$1$};
\node [above] at (11,-4) {$\mu_1$};
\node [above] at (13.5,-4) {$1$};
\node [above] at (16,-4) {$\mu_2$};
\node [above] at (19,-4) {$\cdots$};
\node [above] at (20.5,-4) {$1$};
\node [above] at (23,-4) {$\mu_h$};
\end{tikzpicture}
\end{figure}

In this case, if $|\mu_h| = l \le n-2$, considering for instance the word $u = 1  \mu_h  2^{n-l-2}  0 \in \mathcal{A}_q(n)$ we can clearly see that $1  \mu_h$ is a cross-bifix between $w$ and $u$, and then $W$ is not cross-bifix-free. On the other hand, if $|\mu_h| = n-1$, then necessarily $h = 1$ and $w = 1  \mu_1$. So, $w$ can be written as $w = 1 \alpha \beta$, where $\alpha \in \mathcal{M}_{q-2}(i)$, $\beta \in \hat{\mathcal{M}}_{q-2}(n-i-1)$ with $i > \lfloor\frac{n}{2}\rfloor$ (otherwise $w \in \mathcal{B}_q(n)$). In this case, for instance, the word $\beta  1  2^{i-1}  0 \in \mathcal{A}_q(n)$ has a cross-bifix with $w$, thus $W$ is not a cross-bifix-free-set.\\

If $h < 0$, $w$ can be written as (see Fig. \ref{fig:hneg})
\[ w = \mu_{-h} \; 0 \; \cdots \; \mu_2 \; 0 \; \mu_1 \; 0 \; \phi \]
where $\phi$ is a word satisfying $|\phi|_1 = |\phi|_0$ and ending with $0$, and $\mu_1, \dots, \mu_{-h}$ are $(q-2)$-colored Motzkin words with $\mu_{-h}$ non-empty as $w$ begins with a non-zero letter.

\begin{figure}[h!]
\caption{Graphical representation of $w$, in the case $h < 0$}
\label{fig:hneg}
\centering
\begin{tikzpicture}[scale=0.3]
\draw (0,0) to [out=45, in=180] (2,2) to [out=0, in=135] (4,0) -- (0,0);
\draw (4,0) -- (5,-1);
\draw [dashed] (5,-1) -- (7,-3);
\draw (7,-3) to [out=45, in=180] (9,-1) to [out=0, in=135] (11,-3) -- (7,-3);
\draw (11,-3) -- (12,-4);
\draw (12,-4) to [out=45, in=180] (14,-2) to [out=0, in=135] (16,-4) -- (12,-4);
\draw (16,-4) -- (17,-5);
\draw (17,-5) to [out=45, in=180] (19,-3) to [out=0, in=180] (23,-7) to [out=0, in=225] (25,-5) -- (17,-5);
\draw [dashed] (0,0) -- (25,0) -- (25,-5);
\node [above] at (2,-9) {$\mu_{-h}$};
\node [above] at (4.5,-9) {$0$};
\node [above] at (6,-9) {$\cdots$};
\node [above] at (9,-9) {$\mu_2$};
\node [above] at (11.5,-9) {$0$};
\node [above] at (14,-9) {$\mu_1$};
\node [above] at (16.5,-9) {$0$};
\node [above] at (21,-9) {$\phi$};
\end{tikzpicture}
\end{figure}

In this case, if $|\mu_{-h}| = l \le n-2$, considering for instance the word $u = 1  2^{n-l-2}  \mu_{-h}  0 \in \mathcal{A}_q(n)$ we can clearly see that $\mu_{-h}  0$ is a cross-bifix between $w$ and $u$, and then $W$ is not cross-bifix-free. On the other hand, if $|\mu_{-h}| = n-1$, then necessarily $h = -1$ and $w = \mu_1  0$. So, $w$ can be written as $w = \alpha \beta \delta  0$, where $\beta \in \hat{\mathcal{M}}_{q-2}(j)$ with $j \ge \lceil \frac{n}{2} \rceil$ (otherwise $w \in \mathcal{C}_q(n)$), and $\alpha, \delta$ any two $(q-2)$-colored Motzkin words of the appropriate length. In this case, for instance, the word $2^{n-j-|\alpha|}  \alpha \beta \in \mathcal{A}_q(n)$ has a cross-bifix with $w$, thus $W$ is not a cross-bifix-free-set.\\

Finally, if $h = 0$, the path associated to $w$ can either remain above $x$-axis or fall below it.

In the first case let $i$, with $\lfloor\frac{n}{2}\rfloor \le i < n$, be the last $x$-coordinate of the path intercepting the $x$-axis. Notice that $i$ can not be less than $\lfloor\frac{n}{2}\rfloor$, otherwise $w \in \mathcal{A}_q(n)$. We can write $w = \alpha  \beta$, where $\alpha$ is a non-empty word in $\mathcal{M}_{q-2}(i)$ and $\beta \in \hat{\mathcal{M}}_{q-2}(n-i)$. We now need to take into consideration two different cases: $i = \lfloor\frac{n}{2}\rfloor$ and $i > \lfloor\frac{n}{2}\rfloor$. In the first case $\alpha \in \hat{\mathcal{M}}_{q-2}(\frac{n}{2})$, otherwise $w \in \mathcal{A}_q(n)$, then, for instance, the word $2^{n/2} \alpha \in \mathcal{A}_q(n)$ has a cross-bifix with $w$. In the latter case, for instance, the word $\beta 2^{i-1} 0 \in \mathcal{C}_q(n)$ has a cross-bifix with $w$, so that $W$ is not a cross-bifix-free-set.

In the other case the path associated to $w$ crosses the $x$-axis. Let $i$, with $0 < i < n$, be the first $x$-coordinate of the path crossing $x$-axis. We can write $w = \alpha  0  \phi$, where $\alpha$ is a non-empty word in $\mathcal{M}_{q-2}(i)$. In this case, for instance, the word $1  2^{n-i-2}  \alpha  0 \in \mathcal{A}_q(n)$ has a cross-bifix with $w$, then $W$ is not a cross-bifix-free-set.\\

We proved that $\mathcal{CBFS}_q(n)$ is a non-expandable cross-bifix-free set on $BF_q(n)$, for any $q \ge 3$ and $n \ge 3$.

\cvd

\section{Sizes of Cross-Bifix-Free sets for Small Lengths}

In this section we present some interesting results concerning the size of $\mathcal{CBFS}_q(n)$ compared to the ones in \cite{11}.

For fixed $n$ and $q$, we recall that the size of $q$-ary cross-bifix-free sets given in \cite{11} is obtained by
$$S(n,q) = \max\{(q-1)^2 F_{k,q}(n-k-2): 2 \leq k \leq n-2 \}$$ which is proved to be nearly optimal.

In Table III is shown the values of $S(n,q)$ and $|\mathcal{CBFS}_q(n)|$ for $3\leq q \leq 6$ and $n \leq 16$. For the initial values of $n$, we can observe that the sizes obtained by our construction are greater than the size $S(n,q)$. In particular, the number of the initial values of $n$ for which $|\mathcal{CBFS}_q(n)|$ is greater grows with $q$ and this trend can be easily verified by experimental results.

\begin{table}
\centering
\caption{Comparing the values from \cite{11} with $\mathcal{CBFS}_q(n)$, for $3\leq q \leq 6$}
\scriptsize
\begin{tabular}{c | r r  | r r }
\toprule

$n$ & $|\mathcal{CBFS}_3(n)|$ & $S(n,3)$ & $|\mathcal{CBFS}_4(n)|$ & $S(n,4)$\\
\midrule
3 & 4 & 4 & 9 & 9  \\
4 & \textbf{7} & 4 & \textbf{25} & 9 \\
5 & \textbf{16} & 12 & \textbf{72} & 36  \\
6 & \textbf{36} & 32 & \textbf{223} & 135 \\
7 & 87 & 88 & \textbf{712} & 513 \\
8 & 210 & 240 & \textbf{2 334} & 1 944 \\
9 & 535 & 656 & \textbf{7 868} & 7 371 \\
10 & 1 350 & 1 792 & 26 731 & 27 945 \\
11 & 3 545 & 4 896 & 93 175 & 105 948 \\
12 & 9 205 & 13 376 & 324 520 & 401 679 \\
13 & 24 698 & 36 544  & 1 157 031 & 1 522 881 \\
14 & 65 467 & 99 840 & 4 104 449 & 5 773 680 \\
15 & 178 375 & 272 768 & 14 874 100 & 21 889 683 \\
16 & 480 197 & 745 216 & 53 514 974 & 82 990 089 \\
\bottomrule
\toprule
$n$ & $|\mathcal{CBFS}_5(n)|$ & $S(n,5)$ & $|\mathcal{CBFS}_6(n)|$ & $S(n,6)$\\
\midrule
3 & 16 & 16 & 25 & 25 \\
4 & \textbf{61} & 16 & \textbf{121} & 25 \\
5 & \textbf{224} & 80 & \textbf{550} & 150 \\
6 & \textbf{900} & 384 & \textbf{2 739} & 875 \\
7 & \textbf{3 595} & 1 856 & \textbf{13 260} & 5 125 \\
8 & \textbf{15 014} & 8 960 & \textbf{67 740} & 30 000  \\
9 & \textbf{63 135} & 43 264 & \textbf{342 676} & 175 625 \\
10 & \textbf{271 136} & 208 896 & \textbf{1 787 415} & 1 028 125 \\
11 & \textbf{1 178 677} & 1 008 640 & \textbf{9 324 647} & 6 018 750 \\
12 & \textbf{5 167 953} & 4 870 144 & \textbf{49 456 240} & 35 234 375 \\
13 & 22 986 100 & 23 515 136 & \textbf{263 776 127} & 206 265 625 \\
14 & 102 403 229 & 113 541 120 & \textbf{1 417 981 855} & 1 207 500 000 \\
15 & 463 098 075 & 548 225 024 & \textbf{7 688 015 908} & 7 068 828 125 \\
16 & 2 089 302 415 & 2 647 064 576 & \textbf{41 785 951 916} & 41 381 640 625\\

\bottomrule
\end{tabular}
\end{table}

In order to improve the values of the size $S(n,q)$ for the initial size of $n$,
we can consider the following expression
$$S^*(n,q) = \max\{(q-1)^2 F_{k,q}(n-k-2): 1 \leq k \leq n-2 \},$$
where $k$ can assume also the value 1. When $k=1$, in the case of small $n$ and large $q$, we obtain cross-bifix-free sets
having cardinality greater than the one proposed in \cite{11}.

In Table IV is shown the values of $S^*(n,q)$ and $|\mathcal{CBFS}_q(n)|$ for $3\leq q \leq 6$ and $n \leq 16$.
Also in this situation, we can observe that the sizes obtained by our construction are greater than the size $S(n,q)$ in a range of values of $n$. In particular, the range of values of $n$ for which $|\mathcal{CBFS}_q(n)|$ is greater grows with $q$ and this trend can be easily verified by experimental results.

\begin{table}
\centering
\caption{Comparing the values from $S^*(n,q)$ with $\mathcal{CBFS}_q(n)$, for $3\leq q \leq 6$}
\scriptsize
\begin{tabular}{c | r r  | r r }
\toprule

$n$ & $|\mathcal{CBFS}_3(n)|$ & $S^*(n,3)$ & $|\mathcal{CBFS}_4(n)|$ & $S^*(n,4)$\\
\midrule
3 & 4 & 4 & 9 & 9  \\
4 & 7 & 8 & 25 & 27 \\
5 & 16 & 16 & 72 & 81  \\
6 & \textbf{36} & 32 & 223 & 243 \\
7 & 87 & 88 & 712 & 729 \\
8 & 210 & 240 & \textbf{2 334} & 2 187 \\
9 & 535 & 656 & \textbf{7 868} & 7 371 \\
10 & 1 350 & 1 792 & 26 731 & 27 945 \\
11 & 3 545 & 4 896 & 93 175 & 105 948 \\
12 & 9 205 & 13 376 & 324 520 & 401 679 \\
13 & 24 698 & 36 544  & 1 157 031 & 1 522 881 \\
14 & 65 467 & 99 840 & 4 104 449 & 5 773 680 \\
15 & 178 375 & 272 768 & 14 874 100 & 21 889 683 \\
16 & 480 197 & 745 216 & 53 514 974 & 82 990 089 \\
\bottomrule
\toprule
$n$ & $|\mathcal{CBFS}_5(n)|$ & $S^*(n,5)$ & $|\mathcal{CBFS}_6(n)|$ & $S^*(n,6)$\\
\midrule
3 & 16 & 16 & 25 & 25 \\
4 & 61 & 64 & 121 & 125 \\
5 & 224 & 256 & 550 & 625 \\
6 & 900 & 1 024 & 2 739 & 3 125 \\
7 & 3 595 & 4 096 & 13 260 & 15 625 \\
8 & 15 014 & 16 384 & 67 740 & 78 125  \\
9 & 63 135 & 65 536 & 342 676 & 390 625 \\
10 & \textbf{271 136} & 262 144 & 1 787 415 & 1 953 125 \\
11 & \textbf{1 178 677} & 1 048 576 & 9 324 647 & 9 765 625 \\
12 & \textbf{5 167 953} & 4 870 144 & \textbf{49 456 240} & 48 828 125 \\
13 & 22 986 100 & 23 515 136 & \textbf{263 776 127} & 244 140 625 \\
14 & 102 403 229 & 113 541 120 & \textbf{1 417 981 855} & 1 220 703 125 \\
15 & 463 098 075 & 548 225 024 & \textbf{7 688 015 908} & 7 068 828 125 \\
16 & 2 089 302 415 & 2 647 064 576 & \textbf{41 785 951 916} & 41 381 640 625\\

\bottomrule
\end{tabular}
\end{table}

\newpage
\section{Conclusions and further developments}
In this paper, we introduce a general constructive method for cross-bifix-free sets in the $q$-ary alphabet based
upon the study of lattice paths on the discrete plane. This
approach enables us to obtain the cross-bifix-free set $\mathcal{CBFS}_q(n)$
having greater cardinality than the ones proposed in \cite{11}, for the initial values of $n$.

Moreover, we prove that $\mathcal{CBFS}_q(n)$ is a non-expandable
cross-bifix-free set on $BF_q(n)$, i.e. $\mathcal{CBFS}_q(n) \cup \{w\}$ is
not a cross-bifix-free set on $BF_q(n)$, for any $w \in
BF_q(n) \backslash \mathcal{CBFS}_q(n)$.

The non-expandable property is
obviously a necessary condition to obtain a maximal
cross-bifix-free set on $BF_q(n)$, anyway the problem of determine
maximal cross-bifix-free sets is still open and no general solution has been found yet.

%
%
%

\end{document}